\documentstyle{amsppt}
\overfullrule=0pt

\magnification=\magstep1
 \pagewidth{6.0 true in}
\pageheight{9.0 true in}

\hoffset 0in % change for fine tuning
\voffset -0.10in

\parindent 0pt
\parskip 14pt
\document

\topmatter
\title Remarks on weakly pseudoconvex boundaries
\endtitle
\author Judith Brinkschulte \\
C. Denson Hill \\
Mauro Nacinovich
\endauthor
\address Judith Brinkschulte - Universit\'e de Grenoble 1, Institut
Fourier, B.P. 74,
38402 St. Martin d'H\`eres, France\endaddress
\email brinksch\@ujf-grenoble.fr \endemail
\address C.Denson Hill - Department of Mathematics, SUNY at Stony Brook,
Stony Brook NY 11794, USA \endaddress
\email dhill\@math.sunysb.edu \endemail
\address Mauro Nacinovich - Dipartimento di Matematica -
Universit\`a di Roma "Tor Vergata" - via della Ricerca Scientifica -
00133 - Roma - Italy \endaddress
\email nacinovi\@mat.uniroma2.it\endemail
\abstract In this paper, we consider the boundary $M$ of a weakly pseudoconvex
domain in a Stein manifold. We point out a striking difference between the
local cohomology and the global cohomology of $M$, and illustrate this with an
example. We also discuss the first and second Cousin problems, and the strong
Poincar\'e problem for $CR$ meromorphic functions on the weakly pseudoconvex
boundary $M$.
\endabstract
\subjclass 32V15, 32V05\endsubjclass
\keywords weakly pseudoconvex boundaries, $CR$ meromorphic functions,
tangential Cauchy-Riemann cohomology
\endkeywords
\endtopmatter

\bigskip

Throughout this paper we shall be considering the following situation:
 Let $X$ be a Stein manifold of complex dimension $n+1$, with $n\geq 1$. Let
$\Omega$ be a relatively compact connected domain in $X$, with a smooth
($\Cal{C}^{\infty}$) weakly pseudoconvex boundary $M$. Such an $M$ is then a
connected smooth $CR$ manifold of type $(n,1)$, and we shall refer to it as
being a {\it weakly pseudoconvex boundary}.\smallskip

It is possible to define the tangential $CR$ cohomology groups on $M$, 
both for smooth tangential forms, and for currents, see [AH1], [AH2], [NV].
In a vastly more general situation than for weakly pseudoconvex boundaries, it
was shown in [BHN] that these cohomology groups are either zero, or else must
be infinite dimensional. We discuss below which of these two situations arises
for the special situation of the $M$ under consideration here. \smallskip

In this paper, we point out a striking difference between the local cohomology
and the global cohomology of $M$, and illustrate this with an example.
We also discuss the first and second Cousin problems, and the strong
Poincar\'e problem for $CR$ meromorphic functions on the weakly pseudoconvex
boundary $M$.

\bigskip
\noindent
{\bf \S 1 Local boundary cohomology.}
\smallskip

Let $U$ be an open set on the boundary $M$. We denote by $H^{p,j}(U)$ the 
$j$-th cohomology group of the tangential Cauchy-Riemann operator
$\overline\partial_{M}$ acting on smooth tangential forms of type $(p,*)$ in
$U$, and we denote by $H^{p,j}_{cur}(U)$  the analogous cohomology groups of
$\overline\partial_{M}$ acting on currents of type $(p,*)$ in $U$. Here $0\leq
p\leq n+1$ and $0\leq j\leq n$.

\smallskip

Let $x_0$ be a point on $M$. We consider the localizations, at $x_0$, 
of these cohomology groups: 
$\displaystyle H^{p,j}(\{x_0 \})= \lim_{@>>{U\ni x_0}>} H^{p,j}(U)$, \quad 
$\displaystyle H^{p,j}_{cur}(\{x_0\})= \lim_{@>>{U\ni x_0}>} H^{p,j}_{cur}(U)$.

\medskip

\proclaim{Theorem 1}
Let $x_0$ be a point in $M$ at which the Levi form of $M$ 
has $m$ positive, and $\ell$ zero eigenvalues, with $m+\ell=n$.
Then for all $p$ with $0\leq p\leq n+1$, we have:
\roster
\item"($i$)" 
$H^{p,j}(\{x_0\}) =H^{p,j}_{cur}(\{x_0\})=  0$ for $\ell< j<m$;
\item"($ii$)" 
$H^{p,0}( \{x_0\})$,  $H^{p,0}_{cur}(\{x_0\})$, $H^{p,m}(\{x_0\})$ 
and  $H^{p,m}_{cur}(\{x_0\})$ are infinite dimensional, and the natural map
$$H^{p,m}(\{x_0\}) \longrightarrow H^{p,m}_{cur}(\{x_0\})$$\smallskip
has infinite dimensional image.
\endroster
\endproclaim

\medskip

The vanishing results in (i) were proved for the cohomology for smooth forms 
in [AH2], and these results were extended 
to the cohomology for currents in [NV]. 
In fact in [AH2] and [NV] more general results 
were obtained, which allow also some negative eigenvalues; 
but here we have stated only the special case which pertains 
to weakly pseudoconvex hypersurfaces.\smallskip

 Statement (ii) is obvious when $m=0$ because $\Cal{O}_{x_0}$ is infinite dimensional. When $m=n=1$, it is equivalent to the classical non solvability result of Lewy [L]. For $m=n\geq 1$, it was proved in [AH2] for the smooth case, using a geometrical argument. A completely different proof, which also works in higher codimension, was given in [AFN]. There it was required that the Levi form of $M$ at $x_0$ be nondegenerate in some characteristic conormal direction. The latter argument was recently extended in [HN2] to allow some zero eigenvalues of the Levi form in the same characteristic conormal direction; the statement in (ii) is the special case which pertains to weakly pseudoconvex hypersurfaces.
\medskip
 
\noindent {\bf Interpretation of Theorem 1.}
\smallskip
Let us interprete the results of Theorem 1 in 
terms of solvability or non solvability 
of the system of partial differential equations
$$\overline\partial_M u =f \tag S$$
where the unknown $u$ is a form of type $(p,j-1)$, and the right hand side
$f$ is a form of type $(p,j)$, satisfying the integrability conditions:  
$$\overline\partial_M f = 0.\tag C$$
Note that the above system (S) consists of 
$\left(n+1\atop p\right)+ \left(n\atop j\right)$
first order linear partial differential equations, 
with smooth variable coefficients and no zero order terms, 
which are to be solved for
$\left( n+1\atop p\right) + \left( n\atop j-1\right)$
 unknowns $u$; and the given right hand side $f$ is required to satisfy
$\left(n+1\atop p\right)+ \left(n\atop j+1 \right)$
 similar compatibility conditions. 
Thus the system may be overdetermined, determined, or underdetermined, 
depending on the value of $j$, $1\leq j\leq n$.\smallskip
1. Consider first $H^{p,j}(\{x_0\})=0$: This means that given any open 
neighborhood $U$ of $x_0$ in $M$, 
and any smooth $\overline\partial_M$-closed $(p,j)$-form $f$ in $U$, 
there exist an open neighborhood $V_f$ of $x_0$, 
with $V_f \subset U$, and a smooth $(p,j-1)$-form $u$ 
in $V_f$ such that $\overline\partial_M u=f$ in $V_f$. 
In this situation, it follows via the Baire category theorem and the 
open mapping theorem, that $V_f$ can be chosen independently of $f$; 
i.e., given any open neighborhood $U$ of $x_0$ in $M$, 
there exists an open neighborhood $V$ of $x_0$, 
with $V \subset U$, such that for any smooth $\overline\partial_M$-closed 
$(p,j)$-form $f$ in $U$, there exists a smooth $(p,j-1)$-form $u$ 
in $V$ such that $\overline\partial_M u=f$ in $V$ [HN2].
\smallskip
2. Consider next $H^{p,j}_{cur}(\{x_0\})=0$: The first part of the discussion 
is the same of as above, with smooth forms replaced by currents. 
However in this case the situation concerning the uniform choice of the 
neighborhood $V_f$ is different. Namely, given any neighborhood $U$ of $x_0$ 
in $M$, and any integer $k\geq 0$,
there exist an open neighborhood $V$ of $x_0$ in $M$, with 
$V\subset U$, and an integer $m\geq 0$, such that for any  
$\overline\partial_M$-closed $(p,j)$-current $f$ 
of order $k$ in $U$, there is a 
$(p,j-1)$-current $u$ of order $m$ in $V$, such that 
$\overline\partial_M u=f$ in $V$ in the sense of currents [HN2]. 

\smallskip

We can describe the situations 1. and 2. by saying that the Poincar\'e lemma 
for $\overline\partial_M$ in degree $(p,j)$ is valid at $x_0$.
\smallskip
3. The infinite dimensionality statements in (ii) give local non solvability 
results, which are in the spirit of the scalar Lewy example [L], but which 
are now for overdetermined or underdetermined systems.  In this situation we 
have the following interpretation: There exists a fundamental system $\{ U\}$ 
of open neighborhoods of $x_0$ in $M$, and in each $U$ there are infinitely 
many linearly independent smooth $\overline\partial_M$-closed forms $f$ of 
type $(p,m)$, such that there is no smaller open neighborhood of $x_0$ in $M$ 
in which $f$ can be written as the $\overline\partial_M$ of any current $u$ 
of type $(p,m-1)$.\smallskip
We describe this situation by saying that the Poincar\'e lemma for 
$\overline\partial_M$ in degree $(p,m)$ is not valid at $x_0$.

\bigskip
%%%%%%%%%%%%%%%%%%%%%%%%%%%%%%%%%%%%%%%%%%%%%%%%%%%%%%%%%%%%%%%%%%%%%%%%%%%%%%

\noindent
{\bf \S 2 Global boundary cohomology.}
\smallskip

Next we take $U=M$ and consider the global tangential $CR$ cohomology 
groups on $M$. 

\medskip
\proclaim{Theorem 2}
Assume that $M$ is a weakly pseudoconvex boundary, as in the introduction. 
Then for all $p$ with $0\leq p\leq n+1$ we have:\smallskip
(i) $H^{p,j}(M) = H^{p,j}_{cur}(M) =0$ for $1\leq j\leq n-1$.\smallskip
(ii) $H^{p,0}(M)$, $H^{p,0}_{cur}(M)$, $H^{p,n}(M)$ and $H^{p,n}_{cur}(M)$ 
are infinite dimensional.\smallskip
 If $n\geq 2$, they have a Hausdorff vector space topology.
\endproclaim

This theorem is contained in [B]. Once again, when $j=0$ the statement 
in (ii) is trivial since $X$ is Stein. When $j=n$ the infinite dimensionality 
in (ii) is actually a consequence of the result in Theorem 1. In fact, in 
that case, it was shown in [HN1] that there must exist a point $x_0$ on $M$ 
at which the Levi form of $M$ is positive definite. This gives the infinite 
dimensionality of the local cohomology at $x_0$. However, when $j=n$ we are 
at the end of the $\overline\partial_M$ complex, so there is no compatibility 
condition to 
be satisfied by the right hand side $f$ of (S); 
hence the result globalizes. For the vanishing of the global 
cohomology in (i), the key result proved in [B] was the vanishing of the 
Dolbeault cohomology on $\Omega$ with zero Cauchy data on $M$. The result then
follows by some classical isomorphisms [AH1]. 

\medskip

It follows from Theorem 2 that if at a point $x_0$ in $M$, there is a smooth 
form or current of some type $(p,j)$ with $1\leq j\leq n-1$, which is 
$\overline\partial_M$-closed in some neighborhood of $x_0$, such that one has 
local non solvability for $f$ at $x_0$, then $f$ has no 
$\overline\partial_M$-closed extension to all of $M$.\medskip
%\noindent{\bf Example.}\quad
%Let $z_0,z_1,\hdots,z_n$ be the holomorphic coordinates of $\Bbb C^{n+1}$
%and let
%$$\Omega=\{(z_0,z_1,\hdots,z_n)\in\Bbb C^{n+1}\, | \, 
% |z_0|^{2h}+|z_1|^{2h}+\cdots +|z_n|^{2h}<1\}\, ,$$
%for an integer $m\geq 2$. Then the boundary of $M$ is a 
%weakly pseudoconvex boundary. For each $j=0,1,\hdots,n$ denote by
%$M_j$ consisting of the points $(z_0,z_1,\hdots,z_n)\in M$
%such that $\#\{i\, | \, z_i=0\}=n-j$. Note that $M_j\neq\emptyset$
%for all $j=0,1,\hdots,n$.
%The Levi form of $M$ is
%semidefinite at each point of $M$ and has rank $j$ at points of $M_j$.
%According to [HN2], the Poincar\'e lemma fails in degree $(p,j)$ at
%all points of $M_j$, while according to [B]the global cohomology in
%degree $(p,j)$ vanishes when $1\leq j\leq n-1$.

\medskip
\noindent {\bf Example.} Let $z=(z_0, z_1)$ be coordinates in $\Bbb C^{2}$, 
$w=(w_1,\ldots ,w_{n-1})$ be coordinates in $\Bbb C^{n-1}$. Consider the egg 
in $\Bbb C^{n+1}$ defined by
$$\{ \vert z_0\vert^2 +\vert z_1\vert^2  + \vert w_1\vert^{2m} +\ldots +
\vert w_{n-1}\vert^{2m} < 1\},$$
for an integer $m\geq 2$. It has a weakly pseudoconvex boundary $M$. At each
point $x_0$ on $M$,  there is some value of $j$ with $1\leq j\leq n$ such that
for all  $0\leq p\leq n+1$, the Poincar\'e lemma in degree $(p,j)$ is not
valid  at $x_0$. On the other hand, for any choice of $(p,j)$, 
with $1\leq j\leq n$, there exist points $x_0$ on $M$ at which the
 Poincar\'e lemma in degree $(p,j)$ fails. In fact, for $r=0,1,\ldots ,n-1$, 
let $\Sigma_{n-r}$ be the set of points on $M$ at which exactly $r$ 
components of $w$ are zero. Then $M=\bigcup_{k=1}^n \Sigma_k$, and at each 
point $x_0$ of $\Sigma_k$, the Levi form of $M$ has $k$ positive and $n-k$ 
zero eigenvalues. Hence the Poincar\'e lemma fails at $x_0$ in degree $(p,k)$.
\smallskip
 Nonetheless, for all $(p,j)$ with $0\leq p\leq n+1$ and $1\leq j\leq n-1$ we 
have the global vanishing $H^{p,j}(M)=H^{p,j}_{cur}(M)=0$.\smallskip

A similar situation prevails if we take instead
$$\Omega=\{(z_0,z_1,\hdots,z_n)\in\Bbb C^{n+1}\, | \, 
|z_0|^{m_{0}}+|z_1|^{m_{1}}+\cdots +|z_n|^{m_{n}}<1\}\, ,$$
with $n>1$ and even integers $m_{j}\geq 4$.

%%%%%%%%%%%%%%%%%%%%%%%%%%%%%%%%%%%%%%%%%%%%%%%%%%%%%%%%%%%%%%%%%%%%%%%%%%%%
%

\bigskip
\noindent
{\bf \S 3 Sheaf cohomology.}
\smallskip
In this section we give some
applications of the above results to sheaf cohomology.\smallskip

Let $\Cal O_M$ denote the sheaf of germs of smooth $CR$ functions on $M$ and 
$\Omega^p_M$ denote the sheaf of germs of smooth $CR$ $p$-forms on $M$, thus 
$\Cal O_M =\Omega^0_M$. Likewise we denote by $\Cal O_M^{\prime}$ the sheaf 
of germs of $CR$ distributions on $M$ and by $\Omega^{\prime p}_M$ the sheaf 
of germs of $CR$ $p$-currents on $M$, so 
$\Cal O^{\prime}_M =\Omega^{\prime 0}_M$.
\medskip

\proclaim{Corollary 1} Assume that $M$ is a weakly pseudoconvex boundary, as 
in the introduction, and $n>1$. Then for all $p$ with $0\leq p\leq n+1$, the  
$\check{C}$ech cohomology groups $H^1 (M,\Omega^p_M)$ and  
$H^1 (M,\Omega^{\prime p}_M)$ vanish.
\endproclaim

We obtain the Corollary from Theorem 2 because, by
the abstract de Rham theorem, the maps
$$H^1 (M,\Omega^p_M)\longrightarrow H^{p,1}(M)$$
and
$$H^1 (M,\Omega^{\prime p}_M) \longrightarrow  H^{p,1}_{cur}(M)$$
are injective.
\medskip

The above result can be extended as follows: Let $\Cal Z^{p,j}_M$ denote the 
sheaf of germs of smooth $\overline\partial_M$-closed $(p,j)$-forms on $M$ 
and let $\Cal Z^{p,j}_{M,cur}$ denote the sheaf of germs of 
$\overline\partial_M$-closed $(p,j)$-currents on $M$.\medskip

\proclaim{Corollary 2}
Assume that $M$ is a weakly pseudoconvex boundary, that at each point the 
Levi form of $M$ has at least $m$ positive eigenvalues, and that $2m>n$. Then 
for all $(p,j)$ with $0\leq p\leq n+1$ and $1\leq j < 2m-n$ the $\check{C}$ech
 cohomology groups $H^j(M,\Cal Z^{p,n-m}_M)$ and 
$H^j(M,\Cal Z^{p,n-m}_{M,cur})$ vanish.
\endproclaim

This Corollary is also obtained by using the abstract de Rham 
theorem, because according to Theorem 1, the Poincar\'e lemma is valid at 
each point of $M$ in degree $(p,s)$ for $n-m < s<m$. Hence we have the 
isomorphisms
$$H^j(M,\Cal Z^{p,n-m}_M)\simeq H^{p,n-m+j}(M)$$
$$H^j(M,\Cal Z^{p,n-m}_{M,cur})\simeq  H^{p,n-m+j}_{cur}(M),$$
and therefore again the result follows by Theorem 2.

%%%%%%%%%%%%%%%%%%%%%%%%%%%%%%%%%%%%%%%%%%%%%%%%%%%%%%%%%%%%%%%%%%%%%%%%%%%%%%%%%%%%%%

\bigskip
\noindent
{\bf \S 4 $CR$ meromorphic functions and distributions.}
\smallskip\nopagebreak

In order to consider $CR$ meromorphic functions on $M$, 
we need to assume that the weak unique continuation principle
is valid for $CR$ functions on $M$.
To this aim we shall assume that our weakly pseudoconvex
boundary $M$ is minimal at each point $x_0$. By definition, 
minimality at $x_0$ means that there does not exist any germ of a 
complex $n$-dimensional manifold lying on $M$ and passing through $x_0$. 
It was shown in [DH] that this is equivalent to the nonexistence of any germ 
of a complex $n$-dimensional variety lying on $M$ and passing through $x_0$. 
According to [Tr], [Tu] the minimality condition at $x_0$ implies that the 
germ of a $CR$ function on $M$ at $x_0$ has a local holomorphic extension to 
at least one side of $M$. In our situation the local holomorphic extension 
must be into $\Omega$, because $M$ is weakly pseudoconvex. This yields the 
weak unique continuation property for $CR$ functions on open subsets of $M$; 
i.e., if a $CR$ function $f$ vanishes on some open subset of $M$, then it 
vanishes throughout the connected component of its domain of definition. 
These results are valid for $CR$ distributions $f$, as well as for smooth 
$CR$ functions $f$ on $M$.\medskip
Because of the assumption that $M$ is minimal at each point, we have the 
following result [HN3] concerning the zero locus of $CR$ functions: Let 
$f$, $f\not\equiv 0$, be a continuous $CR$ function in a connected open set 
$U$ on $M$, and set ${\Cal Z_f} =\{x\in U\mid f(x)=0\}$. Then ${\Cal Z_f}$ 
does not disconnect $U$.
\medskip

Let ${\Cal O_M}(U)$ denote the ring of smooth ($\Cal C^{\infty}$) $CR$ 
functions in a not necessarily connected open set $U$; ${\Cal O_M}(U)$ 
is an integral domain if $U$ is connected. Let $\Delta (U)$ be the subset 
of ${\Cal O_M} (U)$ of divisors of zero; i.e., $\Delta (U)$ is the set of 
those $CR$ functions on $U$ which vanish on some open connected component 
of $U$. Let ${\Cal M_M}(U)$ be the quotient ring of ${\Cal O_M}(U)$ with 
respect to ${\Cal O_M} (U) \setminus \Delta (U)$.
 This means that ${\Cal M_M}(U)$ is the set of equivalence classes of 
pairs $(f,g)$ with $f\in{\Cal O_M}(U)$ and 
$g\in{\Cal O_M}(U)\setminus \Delta (U)$. 
The equivalence relation $(f,g) \sim (f^{\prime}, g^{\prime})$ 
is defined by $fg^{\prime} =f^{\prime} g$. If $V\subset U$ is an 
inclusion of open sets, the restriction map 
$r_U^V: {\Cal O_M}(U) \rightarrow {\Cal O_M}(V)$ 
sends ${\Cal O_M}(U)\setminus \Delta (U)$ into 
${\Cal O_M}(V)\setminus \Delta (V)$ and thus induces a homomorphism of rings
$$r_V^U: {\Cal M_M}(U)\longrightarrow {\Cal M_M}(V)$$
We obtain in this way a presheaf of rings. 
We shall call the corresponding sheaf $\Cal M_M$ 
the {\it sheaf of germs of $CR$ meromorphic functions on $M$, in the smooth 
category}. By a $CR$ meromorphic function on $U$ we mean a continuous section 
of $\Cal M_M$ over $U$. In particular, ${\Cal M_M}(M)$ is a field since $M$ 
is connected.
\medskip

Analogously we denote the $CR$ distributions on $U$ by 
${\Cal O_M^{\prime}}(U)$.
 Since $M$ is minimal at each point, we can associate to each open set $U$ in
$M$ a corresponding open set $\widetilde{U}$ in $\overline \Omega$, with
$\widetilde{U}\cap M =U$, such that each $f\in {\Cal O_M^{\prime}}(U)$ has a
unique holomorphic extension $\tilde{f}$ to the interior of $\widetilde{U}$,
having the boundary value $f$ on $U$, in the sense of distributions. If
$f,g\in {\Cal O_M^{\prime}}(U)$ then the product $fg$ may not be defined in
$U$ as a distribution, but it makes sense as a ``hyperfunction'' in $U$,
defined as the ``boundary value'' of $\tilde{f}\tilde{g}$. As before, let
$\Delta^{\prime}(U)$ be the subset of ${\Cal O_M^{\prime}}(U)$ consisting of
those $CR$ distributions in $U$ which vanish on some open connected component
of $U$. We define ${\Cal M_M^{\prime}}(U)$ as the equivalence class of pairs
$(f,g)$, with $f\in {\Cal O_M^{\prime}}(U)$ and $g\in {\Cal
O_M^{\prime}}(U)\setminus \Delta^{\prime}(U)$, with the equivalence relation
$(f_1 ,g_1) \sim (f_2 ,g_2)$ defined by $f_1 g_2 = f_2 g_1$; or what is the
same, by $\tilde{f}_1\tilde{g}_2 = \tilde{f}_2 \tilde{g}_1$. 
To verify that this is an equivalence relation, the transitivity should be
checked on the interior of $\widetilde{U}$ within the field of ordinary
meromorphic functions.\par
If $V\subset U$ is an inclusion of open sets, the restriction map
$r_U^V: {\Cal O_M^{\prime}}(U) \rightarrow {\Cal O_M^{\prime}}(V)$ sends
${\Cal O_M^{\prime}}(U)\setminus \Delta^{\prime} (U)$ into ${\Cal
O_M^{\prime}}(V)\setminus \Delta^{\prime} (V)$ and thus induces a homomorphism
of Abelian groups $$r_V^U: {\Cal M_M^{\prime}}(U)\longrightarrow {\Cal
M_M^{\prime}}(V)$$ We obtain in this way a presheaf of Abelian groups . We
shall call the corresponding sheaf $\Cal M_M^{\prime}$ the {\it sheaf of germs
of $CR$ meromorphic distributions} on $M$. By a $CR$ meromorphic distribution
on $U$ we mean a continuous section of $\Cal M_M^{\prime}$ over $U$. Note
that if $F\in {\Cal M_M^{\prime}}(U)$, and $F\not\equiv 0$ on any connected
component of $U$, then $F$ has a multiplicative inverse $F^{-1}$ in 
${\Cal M_M^{\prime}}(U)$. \medskip

Because of the assumption that $M$ is minimal at each point, we have the 
following 
result [HN3] concerning $CR$ meromorphic distributions: Let 
$g$, $g\not\equiv
0$, be a continuous $CR$ function in a connected open set $U$ on $M$, and set
${\Cal Z_g} = \{ x\in U\mid g(x)=0\}$. Let $F$ be a continuous $CR$ function
defined in $U\setminus {\Cal Z_g}$, which, locally in $U$, satisfies $F =
O(g^{-k})$ for some $k\geq 0$. Then $F$ can be regarded as a $CR$ meromorphic
distribution in $U$, which is locally the quotient of two continuous $CR$
functions.

%%%%%%%%%%%%%%%%%%%%%%%%%%%%%%%%%%%%%%%%%%%%%%%%%%%%%%%%%%%%%%%%%%%%%%%%%%%%%%%%%%%%%%

\bigskip
\noindent
{\bf \S 5 The first boundary Cousin problem.}
\smallskip

We discuss the analogue on $M$ of the first Cousin problem (Mittag-Leffler problem). This means that we want to find a global $CR$ meromorphic function on $M$ which has prescribed `` principal parts''.
\medskip

\proclaim{Theorem 3} Assume that $M$ is a weakly pseudoconvex boundary, 
which is minimal at each point, and $n>1$.
Given an open covering $\{ U_j\}$ of $M$ and principal parts
 $F_j\in {\Cal M_M}(U_j)$ (resp. ${\Cal M_M^{\prime}}(U_j)$), satisfying $F_i
-F_j\in {\Cal O_M}(U_i\cap U_j)$ (resp. $ {\Cal O_M^{\prime}}(U_i\cap U_j)$)
for all $i,j$, there exists $F\in {\Cal M_M}(M)$ (resp. ${\Cal
M_M^{\prime}}(M)$)  such that $F-F_j\in {\Cal O_M}(U_j)$ (resp. ${\Cal
O_M^{\prime}}(U_j)$) on $U_j$. \endproclaim

The argument is the standard one; namely setting $F_i -F_j =h_{ij} 
\in {\Cal O_M}(U_i\cap U_j)$, we have the cocycle condition $h_{ij}+
h_{jk}+h_{ki}=0$ in $U_i\cap U_j\cap U_k$, so $\{ h_{ij}\}$ represents an
element in the first $\check{C}$ech cohomology group of $M$, with respect to
the covering $\{ U_j\}$, having coefficients in the sheaf $\Cal O_M$.
The natural map  $$H^1 (M,\{ U_j\} ,{\Cal O_M})\longrightarrow H^1 (M,{\Cal
O_M})$$ is injective. As $n>1$, the latter group is zero by Corollary 1; hence
the $1$-cocycle is a coboundary, which means that there exist $f_j \in {\Cal
O_M}(U_j)$ such that $h_{ij}=f_i -f_j$ in $U_i\cap U_j$. Setting $F= F_j -f_j$
on $U_j$, we obtain a well defined global $CR$ meromorphic function on $M$
having the desired principal parts. (The argument is the same for $CR$
meromorphic distributions.) \medskip

%%%%%%%%%%%%%%%%%%%%%%%%%%%%%%%%%%%%%%%%%%%%%%%%%%%%%%%%%%%%%%%%%%%%%%%%%%%%%%%

{\bf \S 6 The second boundary Cousin problem.}
\smallskip

Next we discuss the analogue on $M$ of the second Cousin problem (Weierstrass 
problem).
 This means that we want to find a global $CR$ meromorphic function on $M$
which has prescribed ``poles'' and zeros. Let ${\Cal O_M^*}(U)$ denote the
elements of ${\Cal O_M}(U)$ which do not vanish at any point of $U$. \medskip

\proclaim{Theorem 4} Assume that $M$ is a weakly pseudoconvex boundary, 
which is minimal at each point, $n>1$ and $H^2(M,\Bbb Z )=0$.
Given an open covering $\{ U_j\}$ of $M$ and $F_j\in {\Cal M_M}(U_j)$ 
(resp. ${\Cal M_M^{\prime}}(U_j)$),
 satisfying $F_i /F_j\in {\Cal O^*_M}(U_i\cap U_j)$ for all $i,j$, there
exists $F\in {\Cal M_M}(M)$ (resp. ${\Cal
M_M^{\prime}}(M)$) such that $F/F_j\in {\Cal O_M^*}(U_j)$ on $U_j$.
\endproclaim

Let $F_i /F_j =h_{ij} \in {\Cal O_M^*}(U_i\cap U_j)$; now we have the cocycle 
condition $h_{ij}h_{jk}h_{ki}=1$ in $U_i\cap U_j\cap U_k$, so $\{h_{ij}\}$ 
represents an element in the first $\check{C}$ech cohomology group of $M$, 
with respect to the covering $\{ U_j\}$, having coefficients in the sheaf 
$\Cal O_M^*$. 
As above we need to show that $H^1 (M,{\Cal O_M^*})=\{ 1\}$. Here 
$\Cal O_M^*$ is the multiplicative sheaf of germs of smooth never vanishing 
$CR$ functions on $M$. From the short exact sequence of sheaves
$$ 0\longrightarrow \Bbb Z \longrightarrow {\Cal O_M}
\longrightarrow {\Cal O_M^*}\longrightarrow 0$$
we obtain the long exact sequence 
$$\ldots\longrightarrow H^1 (M,{\Cal O_M}) \longrightarrow  
H^1 (M,{\Cal  O_M^*}) 
\longrightarrow H^2 (M, \Bbb Z )\longrightarrow \ldots$$
The desired triviality of $H^1 (M,{\Cal O_M^*})$ follows from Corollary 1 and 
our hypothesis. 

 Hence the $1$-cocycle is a coboundary, which means that there exist 
$f_j \in {\Cal O_M^*}(U_j)$ such that $h_{ij}=f_i /f_j \in U_i\cap U_j$. 
Setting $F= F_j /f_j$ on $U_j$, we obtain a well defined global meromorphic 
function on $M$ having the desired ``poles'' and zeros. Another 
interpretation would be to say that the divisor given by the data in the 
second boundary Cousin problem corresponds to a smooth $CR$ line bundle 
over $M$ which is trivial. 
\medskip

\noindent{\bf Remark.} Suppose in Theorem 4 we prescribe only zeros, but 
no ``poles'' (a positive divisor). This means that $h_j\in {\Cal O_M}(U_j)$ 
for all $j$. Then we obtain a solution $h\in {\Cal O_M}(M)$. But for 
$n\geq 1$ it is well known [AH1] that such a global $CR$ function $h$ 
on $M$ has a smooth extension $\widetilde{h}$ to $\overline\Omega$, which 
is holomorphic in $\Omega$. Thus {\it we are able to construct a global 
holomorphic function in $\Omega$ that has assigned zeros on the boundary of 
$\Omega$}. This holds true for any weakly pseudoconvex boundary:
we don't need in this case  
the assumption of minimality.\medskip

\noindent {\bf Application.}\smallskip
Let $M$ be a weakly pseudoconvex boundary, $n>1$ and $H^2(M, \Bbb Z)=0$.
Consider a smooth compact $CR$ submanifold $S$ in $M$, of type $(n-1,1)$. 
This means that $S$ has real codimension 2 in $M$ and is transversal to the 
Levi distribution on $M$. Assume that $S$ has local smooth $CR$ defining 
functions. {\it Then $S$ has a global smooth defining function, which is $CR$ 
on all of $M$.} Indeed: By hypothesis on $S$ there exists a covering of $S$ 
by open sets $\omega_j$ in $M$, $j=1,\ldots ,N$, such that 
$S\cap\omega_j =\{ x\in \omega_j\mid h_j (x)=0\}$, with 
$h_j\in {\Cal O_M}(\omega_j)$ and $dh_j\not= 0$ on $\omega_j$. Set 
$\omega_0 = M\setminus S$ and $h_0 \equiv 1\in {\Cal O_M}(\omega_0)$. 
We may then assume that $h_i /h_j\in {\Cal O_M^*}(\omega_i\cap\omega_j )$ 
for $i,j=0,1,\ldots ,N$. According to the above remark, there exists a 
global smooth $CR$ function $h$ on $M$ which is a defining function for $S$.

\bigskip

%%%%%%%%%%%%%%%%%%%%%%%%%%%%%%%%%%%%%%%%%%%%%%%%%%%%%%%%%%%%%%%%%%%%%%%%%%%%%%%%%

{\bf \S 7 The strong boundary Poincar\'e problem.}
\smallskip

By the strong Poincar\'e problem on $M$ we mean the following:
 Given a $CR$ meromorphic function on $M$, we want to write it as the
quotient of two global $CR$ functions on $M$, which are coprime at each point.
In order to formulate this question we need to pass to the real analytic
category. Hence in this section, we will assume that $M$ is real analytic. So
now we denote by ${\Cal A_M}(U)$ the ring of real analytic $CR$ functions in
$U$, and by ${\Cal A\Cal M_M}(U)$ the corresponding quotient ring, as before.
In this way we obtain the {\it sheaf of germs ${\Cal A\Cal M_M}$ of $CR$
meromorphic functions on $M$, in the real analytic category}. For a point
$x_0$ on $M$ we denote by $\Cal A_{M,x_0}$ the stalk at $x_0$ in the sheaf
$\Cal A_M$ of germs of real analytic $CR$ functions on $M$. By the
Cauchy-Kowalewski theorem, we have that ${\Cal A_{M,x_0}} \simeq \Cal
O_{x_0}$. Hence  ${\Cal A_{M,x_0}}$ is a unique factorization domain, since $
\Cal O_{x_0}$ is one. Hence the notion of being coprime at $x_0$ makes
sense.\medskip

\proclaim{Theorem 5}
Assume that $M$ is a real analytic weakly pseudoconvex boundary, 
and that $H^2(\Omega,\Bbb Z )=0$. Let $F\in {\Cal A\Cal M_M}(M)$ be a $CR$
meromorphic function on $M$. Then $F=\dsize\frac{G}{H}$, 
where $G,\,H\in {\Cal A_M}(M)$
are real analytic $CR$ functions on $M$, coprime at each point.
\endproclaim

\medskip
By definition $F$ is locally the quotient of real analytic $CR$ functions
 on $M$. It then follows from the Cauchy-Kowalewski theorem that $F$ extends
to a meromorphic function on an open collar neighborhood $V$ of $M$,
that we still denote by the same letter $F$. 
\par
According
to a theorem of [KS], $F$ can be written as a quotient $G_1 /H_1$ with $G_1$
and $H_1$ holomorphic on $V$, but not necessarily coprime. This is because the
envelope of meromorphy of $V$ is Stein, and the weak Poincar\'e problem is
solvable on a Stein manifold. Since $X$ is Stein, $G_1$ and $H_1$ extend 
holomorphically to $V\cup\Omega$.
Hence the mermorphic
function on $V$ extends to a meromorphic function on $V\cup \Omega$, which we
still denote by $F$. By [DF], $\overline\Omega$ has a fundamental system of
open Stein neighborhoods in $X$. Let $\Omega_1$ be a Stein neighborhood of
$\Omega$, chosen sufficiently small, so that $\Omega_1\subset\Omega\cup V$ and
$H^2(\Omega_1 ,\Bbb Z )=0$. It is however a classical result that the strong
Poincar\'e problem has a solution in $\Omega_1$. Thus there exist holomorphic
functions $G$ and $H$ in $\Omega_1$, which are coprime at each point, whose
quotient represents $F$. It suffices to take the restrictions of $G$ and $H$
to $M$.\medskip

\noindent {\bf Remarks.}\smallskip
\noindent 1. Note that the above proof yields more than is stated in the 
theorem; namely, $G$ and $H$ are holomorphic and coprime at each point of
$\overline\Omega$.\smallskip \noindent 2. If we drop the requirement that $G$
and $H$ be coprime at each point, we obtain what is known as the weak
Poincar\'e problem on $M$. This easier problem is always solvable, in the real
analytic category, for a real analytic $M$ which is the boundary of an
$\Omega$ in a Stein manifold, even if we drop the requirements that $M$ be
weakly pseudoconvex, and the topological condition on $\Omega$.

%%%%%%%%%%%%%%%%%%%%%%%%%%%%%%%%%%%%%%%%%%%%%%%%%%%%%%%%%%%%%%%%%%%%%%%%%%%%%%%%%%%%%%


\Refs
\widestnumber\key{ABC}

\ref\key AFN
\by A.
Andreotti, G. Fredricks, M. Nacinovich
\paper On the absence of Poincar\'e
lemma in tangential Cauchy-Riemann
complexes
\jour Ann. Sc. Norm. Sup.
Pisa
\vol 8
\yr 1981
\pages 365-404
\endref


\ref\key AH1
\by A. Andreotti,
C.D. Hill
\paper E.E. Levi convexity and the Hans Lewy problem I: Reduction to 
vanishing theorems
\jour Ann. Sc. Norm. Sup. Pisa
\vol 26
\yr 1972
\pages 325-363
\endref


\ref\key AH2
\by A. Andreotti,
C.D. Hill
\paper E.E. Levi convexity and the Hans Lewy problem II: Vanishing theorems
\jour Ann. Sc. Norm. Sup. Pisa
\vol 26
\yr 1972
\pages 747-806
\endref

\ref\key B
\by J. Brinkschulte
\paper The Cauchy-Riemann equation with support conditions in domains with Levi-degenerate boundaries
\jour Thesis
\endref

\ref\key BHN 
\by J. Brinkschulte, C.D. Hill, M. Nacinovich
\paper Obstructions to generic embeddings
\jour Pr\'epublication de l'Institut Fourier
\vol 538
\yr 2001
\endref

\ref\key DF
\by K. Diederich, J.E. Fornaess
\paper Pseudoconvex domains with real-analytic boundary
\jour Ann. Math.
\vol 107
\yr 1978
\pages 371-384
\endref


\ref\key DH
\by R. Dwilewicz, C.D. Hill
\paper The conormal type function for $CR$ manifolds
\jour to appear in Pub. Math. Debrecen
\endref


\ref \key KS
\by J. Kajiwara, E. Sakai
\paper Generalization of Levi-Oka's theorem concerning meromorphic functions
\jour Nagoya Math. J.
\vol 29
\yr 1967
\pages 75-84
\endref


\ref\key L
\by H. Lewy
\paper An example of a smooth linear partial differential equation without 
solution
\jour Ann. of Math.
\vol 66
\yr 1957
\pages 155-158
\endref

\ref\key HN1
\by C.D. Hill, M. Nacinovich
\paper A necessary condition for
global Stein immersion of compact $CR$
manifolds
\jour Riv. Mat. Univ.
Parma
\vol 5
\yr 1992
\pages 175-182
\endref

\ref\key HN2
\by C.D.Hill, M.Nacinovich
\paper On the failure of the Poincar\'e lemma for
the $\bar\partial_M$-complex
\jour Quaderni sez. Geometria Dip. Matematica Pisa
\vol 1.260.1329
\yr 2001
\pages 1-10
\endref

\ref\key HN3
\by C.D.Hill, M.Nacinovich
\paper The Jacobian theorem for mapping of pseudoconcave $CR$ hypersurfaces
\jour B.U.M.I
\vol (7) 9-A
\yr 1995
\pages 149-155
\endref


\ref\key NV
\by M. Nacinovich, G. Valli
\paper Tangential
Cauchy-Riemann complexes on distributions
\jour Ann. Mat. Pura Appl.
\vol
146
\yr 1987
\pages 123-160
\endref


\ref\key Tr
\by J.M. Trepreau
\paper Sur le prolongement holomorphe des fonctions $CR$ d\'efinies sur une 
hypersurface r\'eelle de classe ${\Cal C^2}$ dans $\Bbb C^n$
\jour Invent. Math.
\vol 83
\yr 1986
\pages 583-592
\endref


\ref\key Tu
\by Tumanov, A.E.
\paper Extension of $CR$ functions into a wedge from a manifold of finite type
\jour Math. Sb. Nov. Ser.
\vol 136
\yr 1988
\pages 128-139
\endref




\endRefs
\enddocument